\documentclass[10pt]{article}
\font\smallit=cmti10

\usepackage{tabularx,tabu}
\usepackage{amssymb,latexsym,amsmath,amsthm,enumerate} 

\makeatletter

\renewcommand\section{\@startsection {section}{1}{\z@}
{-30pt \@plus -1ex \@minus -.2ex}
{2.3ex \@plus.2ex}
{\normalfont\normalsize\bfseries\boldmath}}

\renewcommand\subsection{\@startsection{subsection}{2}{\z@}
{-3.25ex\@plus -1ex \@minus -.2ex}
{1.5ex \@plus .2ex}
{\normalfont\normalsize\bfseries\boldmath}}

\renewcommand{\@seccntformat}[1]{\csname the#1\endcsname. }

\makeatother

\newtheorem{proposition}{Proposition}

\theoremstyle{definition}

\newtheorem{remark}{Remark}

\def\R{\mathcal{R}}
\def\M{\mathcal{M}}

\begin{document}

\begin{center}
\uppercase{
On remoteness functions of 
\; $\MakeLowercase{k}$-NIM with $\MakeLowercase{k}+1$ Piles \\
in Normal and in Mis\`ere Versions}

\vskip 20pt
{\bf Vladimir Gurvich}\\
{\smallit National Research University Higher School of Economics (HSE), Moscow, Russia}\\
{\tt vgurvich@hse.ru}, {\tt vladimir.gurvich@gmail.com}\\
\vskip 10pt
{\bf Vladislav Maximchuk}\\
{\smallit National Research University Higher School of Economics (HSE), Moscow, Russia}\\
{\tt vladislavmaximchuk3495@gmail.com}\\
\vskip 10pt
{\bf Georgy Miheenkov}\\
{\smallit National Research University Higher School of Economics (HSE), Moscow, Russia}\\
{\tt georg2002g@mail.ru}\\
\vskip 10pt
{\bf Mariya Naumova}\\
{\smallit Rutgers Business School, Rutgers University, Piscataway, NJ, United States}\\
{\tt mnaumova@business.rutgers.edu}\\

\end{center}

\centerline{\bf Abstract}
\noindent
Given integer $n$ and $k$ such that $0 < k \leq n$ 
and  $n$  piles of stones, 
two players alternate turns. 
By one move it is allowed 
to choose any  $k$  piles and remove  
exactly one stone from each.
The player who has to move but cannot is the loser.  
in the normal version of the game and 
(s)he is the winner in the mis\`ere version. 
Cases  $k=1$  and $k = n$  are  trivial. 
For $k=2$  the game was solved for  $n \leq 6$. 
For $n \leq 4$  the Sprague-Grundy function was 
efficiently computed 
(for both 
versions). For $n = 5,6$  a polynomial algorithm 
computing P-positions was obtained 
for the normal version. 
\newline 
Then, for the case  $k = n-1$, 
a very simple explicit rule that 
determines the Smith remoteness function 
was found for the normal version of the game: 
the player who has to move keeps a pile 
with the minimum even number of stones; 
if all piles have odd number of stones 
then (s)he keeps a maximum one, 
while the $n-1$ remaining piles are reduced 
by one stone each in accordance with the rules of the game. 
\newline 
Computations show that the same rule 
works efficiently for the mis\`ere version too.  
The exceptions are sparse and are listed in Section 2.  
Denote a position by  $x = (x_1, \dots, x_n)$.  
Due to symmetry, we can assume wlog that  
$x_1 \leq  \ldots \leq x_n$. 
Our computations partition all exceptions 
into the following three families: 
$x_1$ is even, $x_1 = 1$, and  $x_1 \geq 3$ is odd.
In all three cases we suggest explicit formulas 
that cover all found exceptions, but this is not proven. 
\newline 
\noindent{\bf AMS subjects:  91A05, 91A46, 91A68} 

\section{Introduction} 
\label{s0}
We assume that the reader is familiar 
with basic concepts of impartial game theory
(see e.g., \cite{ANW07,BCG01-04} for an introduction) 
and also with the recent papers \cite{GMMV23}, 
where the normal version of 
game NIM$(n,k)$,  the exact slow NIM,   
was analyzed for the case  $n = k+1$.  
Here we consider the mis\`ere version. 

\subsection{Exact Slow NIM} 
\label{ss00}
Game Exact Slow NIM  
was introduced in \cite{GH15} as follows: 
Given two integers  $n$ and $k$  
such that  $0 < k \leq n$  and  
$n$  piles containing  $x_1, \dots, x_n$  stones.  
By one move it is allowed to reduce any  $k$  piles 
by exactly one stone each. 
Two players alternate turns. 
One who has to move but cannot is the loser 
in the normal version of the game and 
(s)he is the winner of its mis\`ere version. 
In \cite{GH15}, this game was denoted  NIM$^1_=(n,k)$.
Here we will simplify this notation to  NIM$(n,k)$.

Game NIM$(n,k)$  is trivial if  $k = 1$  or  $k = n$. 
In the first case it ends after $x_1 + \ldots + x_n$ moves  
and in the second one---after  $\min(x_1, \ldots, x_n)$ 
moves. In both cases nothing depends on players' skills. 
All other cases are more complicated. 

The game was solved for  $k=2$  and  $n \leq 6$.
In \cite{GHHC20}, an explicit formula 
for the Sprague-Grundy 
(SG) function was found for  $n \leq 4$, 
for both the normal and mis\`ere versions. 
This formula allows us to compute  the SG function in linear time. 
Then, in \cite{CGKPV21} the P-positions 
of the normal version were found for  $n \leq 6$.  
For the subgame with even  $x_1 + \ldots + x_n$ 
a simple formula for the P-positions was obtained, 
which allows to verify in linear time 
if  $x$ is a P-position and, if not, 
to find a move from it to a P-position.  
The subgame with odd  $x_1 + \ldots + x_n$ is  more difficult. 
Still a (more sophisticated)  formula 
for the P-positions was found,  
providing a linear time recognition algorithm.

Further generalizations of exact slow NIM 
were considered in \cite{GMN23,GN23,GN23a}.

\subsection{Case $n = k+1$, the normal version}
\label{ss01}
In \cite{GMMV23} the normal version 
was solved in case $n = k+1$  by the following simple rule: 
\begin{itemize}
\item{(o)} if all piles are odd, 
keep a largest one and reduce all other;
\item{(e)} if there exist even piles, 
keep a smallest one of them and reduce all other.
\end{itemize}

This rule is well-defined and it uniquely determines 
a move in every position $x$\footnote{Obviously, permuting the piles 
with the same number of stones we keep the game unchanged.}. 
The rule and the corresponding moves 
are called the {\em M-rule} and {\em M-moves};  
the sequence of successive M-moves 
is called the {M-sequence}.

Obviously,  $n > 1$ is required. 
If  $n=1$  then  $x_1$  will reach 
an even value in at most one M-move and stop. 
Since this case is trivial, 
we can assume that  $n > 1$  
without any loss of generality (wlog). 

It is also easily seen that 
no M-move can result in a position whose all entries are odd.  
Hence, for an M-sequence,  part (o) of the M-rule 
can be applied at most once, at the beginning; 
after this only part (e) works. 

Given a position  $x = (x_1, \ldots, x_n)$, 
assume that both players follow the M-rule 
and  denote by $\M(x)$  the number of moves 
from  $x$  to a terminal position. 
In \cite{GMMV23} it was proven that  $\M = \R$, where 
$\R$  is the classical {\em remoteness function} 
introduced by Smith \cite{Smi66}. 
Thus, M-rule solves the game and, moreover, 
it allows to win as fast as possible in an N-position and 
to resist as long as possible in a P-position.

A polynomial algorithm computing  $\M = \R$  
(and in particular, the P-positions) is given, 
even if $n$  is a part of the input and 
integers are presented in binary form.

Let us also note that 
an explicit formula for the P-positions is known 
only for $n \leq 4$  
and already for $n=3$ 
it is pretty complicated 
\cite{GHHC20}\cite[Appendix]{GMMV23}.

\subsection{Related versions of NIM} 
By definition, the present game NIM$(n,k)$ 
is the exact slow version  
of the famous Moore's NIM$_k$ \cite{Moo910}.  
In the latter game a player, by one move, 
reduces arbitrarily 
(not necessarily by one stone) 
at most $k$ piles from $n$. 

The case  $k=1$  corresponds to the classical NIM  
whose P-position was found by Bouton 
\cite{Bou901}  for both the normal and mis\`ere versions. 

\begin{remark} 
Actually, the Sprague-Grundy (SG) values
of NIM were also computed in Bouton's paper, 
although were not defined explicitly in general. This 
was done later by Sprague \cite{Spr36} 
and Grundy \cite{Gru39} for 
arbitrary disjunctive compounds 
of impartial games; see also \cite{Smi66,Con76}. 

In fact, the concept of a P-position 
was also introduced by Bouton in \cite{Bou901}, 
but only for the (acyclic) digraph of NIM, 
not for all impartial games. 
In its turn, this is a special case 
of the concept of a kernel, 
which was introduced for arbitrary digraphs 
by von Neumann and Morgenstern \cite{NM44}. 

Also the mis\`ere version 
was introduced by Bouton in \cite{Bou901},  
but only for NIM, not for all impartial games. 
The latter was done by Grundy and Smith 
\cite{GS56}; see also \cite{Con76,Smi66,Gur07,Gur07a,GH18}.
\end{remark} 

Moore \cite{Moo910} obtained an elegant explicit formula 
for the P-positions of NIM$_k$ 
generalizing the Bouton's case  $k=1$. 
Even more generally, the positions 
of the SG-values 0 and 1  were efficiently characterized by 
Jenkins and Mayberry \cite{JM80}; see also 
\cite[Section 4]{BGHMM17}.  

Also in \cite{JM80};
the SG function of NIM$_k$  
was computed explicitly for the case  $n = k+1$ 
(in addition to the case $k=1$).
In general, no explicit formula, 
nor even a polynomial algorithm, 
computing the SG-values 
(larger that 1)  is known. 
The smallest open case: 
2-values for  $n=4$ and $k=2$. 

\medskip 

The remoteness function of  $k$-NIM was recently studied 
in \cite{BGMV23}. 

\medskip 

Let us also mention 
the exact (but not slow) game NIM$^=(n,k)$  \cite{BGHM15,BGHMM17} 
in which exactly $k$ from $n$ piles are reduced 
(by an arbitrary number of stones) in a move. 
The SG-function was efficiently computed 
in \cite{BGHM15,BGHMM17} for $n \leq 2k$. 
Otherwise, even a polynomial algorithm 
looking for the P-positions is not known
(unless $k=1$, of course). 
The smallest open case is  $n=5$  and  $k=2$.

\section{Case $n = k+1$, mis\`ere version}
\label{s1}
Computations show that the same M-rule 
works pretty efficiently also for the mis\`ere version 
of the considered game NIM$(n,k)$  with $n = k+1$. 
A position $x = (x_1, \dots, x_n)$  is called 
an {\em exception} if 
$\R(x) - \R(x') \neq 1$  for the M-move $x \to x'$.
In fact, $\R(x)$ is odd and $\R(x) - \R(x')$ takes values 0 or 2 
for all known exceptions.\footnote{$\R(x)$ and  $\R(x')$ 
can be equal, but cannot be both even.}
The exceptions are sparse and 
satisfy a regular pattern based on 
two parameters: $n$  and  $\min(x_1, \dots, x_n)$.  
However, the complete description of this pattern 
is open; also no proofs are given in this section. 

\subsection{Monotonicity for the entries of positions} 
Recall that  $x_1 \leq \ldots \leq x_n$  is assumed 
for any position  $x = (x_1, \ldots, x_n)$\footnote{We order the entries $x_i$ 
just for convenience;
their permutations do not change the game.}. 
 
However, even if this monotonicity  holds for  $x$, 
it may fail for $x'$, after a move $x \to x'$. 
In this case, we have to restore it 
by permuting entries of  $x'$.

Alternatively, we can make the M-rule slightly stricter, as follows. 
Given a position $x = (x_1, \ldots, x_n)$ 
for which the M-rule is "ambiguous``, 
that is, $x$ contains 

\medskip

(o) several  smallest even entries, or 

\smallskip

(e) several largest odd entries, provided all $x_i$  are odd. 

\medskip

In both cases, among these equal entries 
keep one with the largest index reducing all others by 1.
We will call such M-move 
(as well as the corresponding M-sequence and M-rule) 
{\em strict}.
It is easily seen that a strict M-move $x \to x'$ 
respects the non-decreasing monotonicity of the entries, that is,  
$(x'_1 \leq  \ldots \leq x'_n)$  whenever 
$(x_1 \leq  \ldots \leq x_n)$.
In contrast, every non-strict M-move breaks this monotonicity.   

\subsection{Monotonicity of exceptions}
\label{ss11a} 
A position is called an {\em exception} 
if the M-move is not optimal in it.  
An optimal move in such position is called 
{\em exceptional}. 

\begin{proposition} 
Given integer  $m \geq n$, a position    
$x' = (x_1, \ldots, x_n, \dots x_m)$  
is an exception  whenever   
$x = (x_1, \ldots, x_n)$  is.  
Of course, not vice versa. 
Moreover, the exceptional moves
coincide  in  $x$  and  $x'$.\footnote{More precisely, the entry 
$x_i$  that is kept unchanged by an optimal move 
is the same for $x$ and $x'$; 
furthermore,  $1 \leq i \leq n$.}
\end{proposition} 


Listing exceptional position below,   
by default, we do not include  $x'$  
if  $x$  is already listed;  
in other words, we include only the minimal exceptions. 

\subsection{General properties of exceptions}
\label{ss11b} 
All found minimal monotone exception 
share the following properties: 

\begin{itemize}
\item{[$x_n$-monotone]}  $x_n > x_{n-1}$;  
if position  $x^i = (x_1, \dots, x_{n-1}, x_{n-1} + i)$  
is an exception for some  $i > 0$  
then  $x^i$  is an exception for each  $i > 0$. 

\item{$x_n - x_{n-1} = 1$} in every minimal exception. 

\item{[$x_{n-1}-determining]$} 
$\R(x) = f(x_{n-1}) + 1$, where 
$f(\ell) = 2 \lceil \ell/2 \rceil$, that is, 
$f(\ell) = \ell$  if $\ell$ is even 
and  $f(\ell) = \ell+1$  if  $\ell$  is odd, 
for all integer $\ell \geq 0$. 
\newline 
Thus, $\R(x)$  (and $\R(x')$) are odd 
(and, hence, the first player wins) in every exception.
However, the M-move is losing. 
(It could win but require a larger number of moves.
Yet, such case is not realized in any found exception.)  

\item 
In every exception $x$, 
the optimal move keeps the entry $x_n$ 
if $x_{n-1}$  is even and keeps $x_{n-1}$ 
if it is odd. 
In contrast, the strict M-move, 
vice versa, keeps  $x_{n-1}$ 
if it is even and keeps $x_n$ 
if  $x_{n-1}$  is odd. 

\item{[$0 \leq \R(x) - \R(x') \leq 2$]} 
By definition of the remoteness function, 
we have  $\R(x)-\R(x')=1$  
for each optimal move $x \to x'$  in every impartial game,  
in particular, for every M-move 
in the normal version of NIM$(n,n-1)$. 
In contrast, for its mis\`ere version, 
in every found {\em minimal} exception 
$\R(x) - \R(x')$  takes only 
values 0, when $x_{n-1}$ is even, or 2, when $x_{n-1}$ is odd. 

\bigskip


\end{itemize}

\subsection{Even  $x_1$}
\label{ss12}  
Given  $x_1 = 2i$, 
a position $x = (x_1, \ldots, x_n, \ldots x_m)$  
is an exception if and only if 
$$2i = x_1 = \ldots = x_{i+2} < x_{i+3} = x_n \leq \ldots \leq x_m \; 
\text{where} \;  1 \leq i \leq n-3  \; \text{and} \; n \geq 4.$$ 

Since $x_{n-1}$ is always even, 
the optimal move keeps  $x_n$, while 
the M-move keeps  $x_j$ for any fixed $j < n$. 
Since  $\R(x) = 2i+1$  is odd, the first player always wins. 
Note that both properties agree with Subsection \ref{ss11b}. 

Examples for $x_1 = 2, 4, 6, 8$  are given below 

\begin{tabu} to \textwidth{X[4.4] X X[1.1] X[0.5]}
$2 = x_1 = x_2 = x_3 < x_4 = x_n \leq \dots \leq x_m,$ & $\R(x) = 3,$ & $n=4$;&\\
$4 = x_1 = \ldots = x_4 < x_5 = x_n \leq \dots \leq x_m,$ & $\R(x) = 5,$ & $n=5$;&\\
$6 = x_1 = \ldots = x_5 < x_6 = x_n \leq \dots \leq x_m,$ &  $\R(x) = 7,$ & $n=6$;&\\
$8 = x_1 = \ldots = x_6 < x_7 = x_n \leq \dots \leq x_m,$ &  $\R(x) = 9,$ & $n=7, \text{ etc.}$&\\
\end{tabu}

\vspace{5mm}

\subsection{$x_1 = 1$}
\label{ss13}  
A position $x = (x_1, \ldots, x_n, \ldots, x_m)$  with  
$x_1 = 1$ 
is an exception if and only if 
$$1 = x_1 \leq x_2 < x_3 \leq \ldots \leq x_m,  \; 
\text{here} \; 3 = n \leq m.$$ 

Furthermore, if  $x_2$  is even then 
$\R(x) = x_2 + 1$  and the only 
optimal move keeps  $x_3$, 
while the M-move keeps  $x_2$; 
if  $x_2$  is odd then 
$\R(x) = x_2 + 2$  and, in contrast,  
the only optimal move keeps  $x_2$, 
while the M-move keeps  $x_i$,  
for some $i > 2$.
In both cases  $\R(x)$  is odd and,  
hence, the first player always wins.  
It is easily seen that all these properties 
agree with Section \ref{ss11b}. 

Wlog, we could restrict ourselves by  $n=3$. 
All exceptions with larger  $n$  are implied  by monotonicity. 

\medskip 

Thus, it remains to consider odd  $x_1 \geq 3$. 
Our computer analysis includes only  $x_1 = 3,5,7,9,11,13,15,17$. 
In each case we observe some pattern, 
yet, its extension to arbitrary odd  $x_1$ 
remains an open problem. 

\subsection{Odd $x_1 \geq 5$ with $n=4$}
\label{ss14}  
For any odd $x_1 \geq 5$  
fix an integer  $i \geq 0$  to obtain the following 
two exceptions  $x = (x_1, x_2, x_3, x_4)$:   

$$x_1, \; x_2 = x_1 + 2i, \; x_3 = 2(x_1 + i - 2), \; x_4 = 2x_1 + 2i - 3 = x_3 + 1;$$
$$x_1, \; x_2 = x_1 + 2i, \; x_3 = 2 x_1 + 2i -3,  \; x_4 = 2(x_1 + i - 1) = x_3 + 1.$$

In the first case, $\R(x) = x_4$  and 
the unique optimal move keeps $x_4$, 
while the unique M-move keeps $x_3$; 
in contrast, in the second case, $\R(x) = x_4+1$  and 
the unique optimal move keeps $x_3$, 
while the unique M-move keeps $x_4$. 

Furthermore, the remoteness function is given by formula 
$\R(x) = f(x_3 + 1)$, where 
$f(m) = 2 \lfloor m/2 \rfloor + 1$, that is, 
$f(m) = m$  if $m$ is odd and  $f(m) = m+1$  if  $m$  is even, 
for all $m \geq 0$. 

It is easily seen that  $\R(x)$  is odd for each  $i$; 
hence, the first player always wins. 

\medskip 

Examples for $x_1 = 5,7,9,11$  are given below. 
Notation  $y+$  means ``any number which is greater than or equal to $y$".

\begin{center}
\footnotesize
\begin{tabular}{|ll|ll|ll|ll|}
\hline
$x_1 = 5$ & $\R$ & $x_1 = 7$ & $\R$ & $x_1 = 9$ & $\R$ & $x_1 = 11$ & $\R$\\
\hline
(5,\;5,\;\;6,\;\;7+) & 7& (7,\;7,\; 10, 11+) & 11 & (9,\,\,\;9,\;\,\, 14, 15+) & 15 & (11, 11,\,\;18, 19+) & 19\\
(5,\;5,\;\;7,\;\;8+)& 9& (7,\;7,\; 11, 12+)& 13 & (9,\,\,\;9,\;\,\, 15, 16+) & 17 & (11, 11,\,\;19, 20+) & 21\\
(5,\;7,\;\;8,\;\;9+)& 9& (7,\;9,\; 12, 13+)& 13 & (9, 11,\,\,\;16, 17+) & 17 & (11, 13,\,\;20, 21+) & 21\\
(5,\;7,\;\;9,\;10+)& 11& (7,\;9,\; 13, 14+)& 15 & (9, 11,\,\,\;17, 18+) & 19 & (11, 13,\,\;21, 22+) & 23\\
(5,\;9,\;10, 11+)& 11& (7, 11,\;14, 15+) & 15 & (9, 13,\,\,\;18, 19+) & 19 & (11, 15,\,\;22, 23+) & 23\\
(5, 9,\;11, 12+)& 13& (7, 11,\;15, 16+) & 17 & (9, 13,\,\,\;19, 20+) & 21 &  (11, 15,\,\;23, 24+) &  25\\
\hline
\end{tabular}
\end{center}

\bigskip 

By monotonicity, any such exception  
$x = (x_1, x_2, x_3, x_4)$ 
can be extended to the exceptions $x' = (x'_1, \ldots,  x'_m)$, 
with  $m \geq 5$  and 
$x_i = x'_i$  for $i \leq 4$, while 
$x'_5, \ldots x'_m$  can be chosen arbitrary such that 
$x_4 \leq x'_5 \leq \ldots \leq x'_m$. 
Note also that case  $x_1 = 3$  is considered in Section \ref{ss-15}.

\subsection{Odd $x_1 \geq 7$ with $n=5$}
\label{ss14-5}  

The following families of exceptions were found: 

\bigskip 

\begin{center}
\footnotesize
\begin{tabular}{|l|l|l|l|}
\hline
$x_1 = 7$ & $x_1 = 9$ & $x_1 = 11$ & $x_1 = 13$\\
\hline
(7, 7, \, \;8, \;8, \; 9) & (9,\, 9, \,\; 12, 12, 13) & (11, 11, \, 16, 16, 17) & (13, 13, \, 20, 20, 21)\\
(7, 7, \, \;9, \;9, \,10)& (9,\, 9, \,\; 13, 13, 14) & (11, 11, \, 17, 17, 18) & (13, 13, \, 21, 21, 22)\\
(7, 9, \, 10, 10, 11) & (9, 11, \, 14, 14, 15) & (11, 13, \, 18, 18, 19) &(13, 15, \, 22, 22, 23)\\
(7, 9, \, 11, 11, 12) & (9, 11, \, 15, 15, 16) & (11, 11, \, 19, 19, 20) &(13, 15, \, 23, 23, 24)\\  
(7, 11, 12, 12, 13)& (9, 13, \, 16, 16, 17) & (11, 15, \, 20, 20, 21) &(13, 17, \, 24, 24, 25) \\
\hline
\end{tabular}
\end{center}

\bigskip 
\begin{center}
\footnotesize
\begin{tabular}{|l|l|l|l|}
\hline
$x_1 = 9$ & $x_1 = 13$ & $x_1 = 17$ & $x_1 = 21$\\
\hline
(9, \, 9, \, 9, \, 10, 11+) & (13, 13, 13, \, 16, 17+)  & (17, 17, 17, \, 22, 23+) & (21, 21, 21, \, 28, 29+)\\
(9, \, 9, \, 9, \, 11, 12+) & (13, 13, 13, \, 17, 18+) & (17, 17, 17, \, 23, 24+) & (21, 21, 21, \, 29, 30+)\\
(9, 11, 11, \, 12, 13+) & (13, 15, 15, \, 18, 19+) & (17, 19, 19, \, 24, 25+) & (21, 23, 23, \, 30, 31+)\\
(9, 11, 11, \, 13, 14+) & (13, 15, 15, \, 19, 20+) & (17, 19, 19, \, 25, 26+) & (21, 23, 23, \, 31, 32+)\\
(9, 13, 13, \, 14, 15+) & (13, 17, 17, \, 20, 21+) & (17, 21, 21, \, 26, 27+) & (21, 25, 25, \, 32, 33+)\\
\hline
\end{tabular}
\end{center}

\bigskip 

\begin{center}
\footnotesize
\begin{tabular}{|l|l|l|l|}
\hline
$x_1 = 11$ & $x_1 = 13$ & $x_1 = 15$ & $x_1 = 17$\\
\hline
(11, 11, 13, \, 14, 15+) & (13, 13, 17, \, 18, 19+) & (15, 15, 21, \, 22, 23+)  & (17, 17, 25, \, 26, 27+) \\
(11, 11, 13, \, 15, 16+) & (13, 13, 17, \, 19, 20+) & (15, 15, 21, \, 23, 24+)  & (17, 17, 25, \, 27, 28+) \\
(11, 13, 15, \, 16, 17+) & (13, 15, 19, \, 20, 21+)  & (15, 17, 23, \, 24, 25+)& (17, 19, 27, \, 28, 29+)\\
(11, 13, 15, \, 17, 18+) & (13, 15, 19, \, 21, 22+) & (15, 17, 23, \, 25, 26+)& (17, 19, 27, \, 29, 30+) \\
(11, 15, 17, \, 18, 19+) & (13, 17, 21, \, 22, 23+)& (15, 19, 25, \, 26, 27+)  & (17, 21, 29, \, 30, 31+)\\
\hline
\end{tabular}
\end{center}

\bigskip 

\begin{center}
\footnotesize
\begin{tabular}{|l|l|l|l|}
\hline
$x_1 = 15$ & $x_1 = 17$ & $x_1 = 19$ & $x_1 = 21$\\
\hline
(15, 15, 17, \, 20, 21+)&(17, 17, 21, \,  24, 25+) &(19, 19, 25, \, 28, 29+) &(21, 21, 29, \, 32, 33+)\\
(15, 15, 17, \, 21, 22+)&(17, 17, 21, \,  25, 26+) &(19, 19, 25, \, 29, 30+) &(21, 21, 29, \, 33, 34+)\\
(15, 17, 19, \, 22, 23+)&(17, 19, 23, \,  26, 27+) &(19, 21, 27, \, 30, 31+) &(21, 23, 31, \, 34, 35+) \\
(15, 17, 19, \, 23, 24+)&(17, 19, 23, \,  27, 28+) &(19, 21, 27, \, 31, 32+) &(21, 23, 31, \, 35, 36+) \\
(15, 19, 21, \, 24, 25+)&(17, 21, 25, \,  28, 29+) &(19, 23, 29, \, 32, 33+) &(21, 25, 33, \, 36, 37+)\\
\hline
\end{tabular}
\end{center}

\subsection{Odd $x_1 \geq 9$ with $n=6$}
\label{ss14-6}
The following families of exceptions were found: 

\begin{center}
\footnotesize
\begin{tabular}{|l|l|l|l|}
\hline
$x_1 = 9$ & $x_1 = 11$ & $x_1 = 13$ & $x_1 = 15$\\
\hline
(9, 9, \,\;\, 10, 10, 10, 11+)& (11, 11, \, 14, 14, 14, 15+) & (13, 13, \, 18, 18, 18, 19+) & (15, 15, \, 22, 22, 22, 23+)\\ 
(9, 9,\, \;\, 11, 11, 11, 12+) & (11, 11, \, 15, 15, 15, 16+)& (13, 13, \, 19, 19, 19, 20+) & (15, 15, \, 23, 23, 23, 24+)\\
(9, 11, \, 12, 12, 12, 13+)& (11, 13, \, 16, 16, 16, 17+)& (13, 15, \, 20, 20, 20, 21+)& (15, 17, \, 24, 24, 24, 25+) \\
(9, 11, \, 13, 13, 13, 14+)& (11, 13, \, 17, 17, 17, 18+) & (13, 15, \, 21, 21, 21, 22+)& (15, 17, \, 25, 25, 25, 26+) \\
(9, 13, \, 14, 14, 14, 15+)& (11, 15, \, 18, 18, 18, 19+)& (13, 17, \, 22, 22, 22, 23+)& (15, 19, \, 26, 26, 26, 27+) \\
\hline
\end{tabular}
\end{center}

\begin{center}
\footnotesize
\begin{tabular}{|l|l|l|l|}
\hline
$x_1 = 11$ & $x_1 = 15$ & $x_1 = 19$ & $x_1 = 23$\\
\hline
(11, 11, 11, \, 12, 12, 13+)& (15, 15, 15, \, 18, 18, 19+) & (19, 19, 19, \, 24, 24, 25+)& (23, 23, 23, \, 30, 30, 31+)  \\
(11, 11, 11, \, 13, 13, 14+) & (15, 15, 15, \, 19, 19, 20+)& (19, 19, 19, \, 25, 25, 26+)& (23, 23, 23, \, 31, 31, 32+) \\
(11, 13, 13, \, 14, 14, 15+) & (15, 17, 17, \, 20, 20, 21+) & (19, 21, 21, \, 26, 27, 27+)& (23, 25, 25, \, 32, 32, 33+)  \\
(11, 13, 13, \, 15, 16, 16+)& (15, 17, 17, \, 21, 21, 22+)& (19, 21, 21, \, 27, 27, 28+)& (23, 25, 25, \, 33, 33, 34+) \\
(11, 15, 15, \, 16, 16, 17+)& (15, 19, 19, \, 22, 22, 22+)&  (19, 23, 23, \, 28, 29, 29+)& (23, 27, 27, \, 34, 34, 35+) \\
\hline
\end{tabular}
\end{center}

\begin{center}
\footnotesize
\begin{tabular}{|l|l|l|}
\hline
$x_1 = 15$ & $x_1 = 17$ & $x_1 = 19$\\
\hline
(15, 15, 17, 17, \, 18, 19+)&(17, 17, 21, 21, \, 22, 23+)&(19, 19, 25, 25, \, 26, 27+)\\
(15, 15, 17, 17, \, 19, 20+)&(17, 17, 21, 21, \, 23, 24+)&(19, 19, 25, 25, \, 27, 28+)\\
(15, 17, 19, 19, \, 20, 21+)&(17, 19, 23, 23, \, 24, 25+)&(19, 21, 27, 27, \, 28, 29+)\\
(15, 17, 19, 19, \, 21, 22+)&(17, 19, 23, 23, \, 25, 26+)&(19, 21, 27, 27, \, 29, 30+)\\
(15, 19, 21, 21, \, 22, 23+)&(17, 21, 25, 25, \, 26, 27+)&(19, 23, 29, 29, \, 30, 31+)\\
\hline
\end{tabular}
\end{center}

\begin{center}
\footnotesize
\begin{tabular}{|l|l|l|l|}
\hline
$x_1 = 13$ & $x_1 = 15$ & $x_1 = 17$ & $x_1 = 19$\\
\hline
(13, 13, 15, \, 16, 16, 17+)&(15, 15, 19, \, 20, 20, 21+)&(17, 17, 23, \, 24, 24, 25+) &(19, 19, 27, \, 28, 28, 29+)\\
(13, 13, 15, \, 17, 17, 18+)&(15, 15, 19, \, 21, 21, 22+)&(17, 17, 23, \, 25, 25, 26+) &(19, 19, 27, \, 29, 29, 30+)\\
(13, 15, 17, \, 18, 18, 19+)&(15, 17, 21, \, 22, 22, 23+)&(17, 19, 25, \, 26, 26, 27+)&(19, 21, 29, \, 30, 30, 31+)\\
(13, 15, 17, \, 19, 19, 20+)&(15, 17, 21, \, 23, 23, 24+)&(17, 19, 25, \, 27, 28, 28+)&(19, 21, 29, \, 31, 31, 32+)\\
(13, 17, 18, \, 20, 20, 21+)&(15, 19, 21, \, 24, 24, 25+)&(17, 21, 27, \, 28, 28, 29+) &(19, 23, 31, \, 32, 32, 33+)\\
\hline
\end{tabular}
\end{center}

\begin{center}
\footnotesize
\begin{tabular}{|l|l|l|}
\hline
$x_1 = 13$ & $x_1 = 19$ & $x_1 = 25$\\
\hline
(13, 13, 13, 13, \, 14, 15+)&(19, 19, 19, 19, \, 22, 23+)&(25, 25, 25, 25, \, 30, 31+)\\ 
(13, 13, 13, 13, \, 15, 16+)&(19, 19, 19, 19, \, 23, 24+)&(25, 25, 25, 25, \, 31, 32+)\\
(13, 15, 15, 15, \, 16, 17+)&(19, 21, 21, 21, \, 24, 25+)&(25, 27, 27, 27, \, 32, 33+)\\
(13, 15, 15, 15, \, 17, 18+)&(19, 21, 21, 21, \, 25, 26+)&(25, 27, 27, 27, \, 33, 34+) \\
(13, 17, 17, 17, \, 18, 19+) &(19, 23, 23, 23, \, 26, 27+)&(25, 29, 29, 29, \, 34, 35+) \\  
\hline
\end{tabular}
\end{center}

Note that  the second and last families 
are defined only for  
$x_1 = 4m + 3$  and  $x_1 = 6m + 1$,
respectively, 
where  $m \geq 3$.  

\subsection{Odd  $x_1 \geq 3$  with  $n = \frac{1}{2}(x_1+1) + 2$}  
\label{ss-15}  
For $n - \frac{1}{2}(x_1+1) > 2$ no exceptions were found, 
while in the considered case 
the exceptions are as follows. 
Given an integer $i \geq 0$, 
a position $x = (x_1, \ldots, x_n)$  is an exception if and only if 
$$x_1 + i = x_2 = \ldots = x_{n-1} < x_n \leq \dots \leq x_m, \; i = 0,1, \ldots$$  

Furthermore,  the remoteness function is given by formula 
$\R(x) = f(x_1 + i + 1)$, where function  $f$  was defined 
in Section \ref{ss11b}.

Note that  $x_2 = \ldots = x_{n-1}$  
and this number is even if and only if  $i$  is odd.

Examples for $x_1 = 3, 5, 7, 9, 11$  and, respectively, 
$n = 4, 5, 6, 7, 8$  are given below. 

\begin{center}
\footnotesize
\begin{tabular}{|ll|ll|ll|}
\hline
$x_1 = 3, n = 4$ & $\R$ & $x_1 = 5, n = 5$ & $\R$ & $x_1 = 7, n = 6$ & $\R$\\
\hline
(3,\;\,3, 3, 4+) & 5 & (5,\;5, 5, 5,\;6+)& 7 &(7,\;\;7, \;7,\;7, \;7, \,8+)& 9\\
(3,\;\,4, 4, 5+) & 5 & (5,\;6, 6, 6,\;7+)& 7 &(7,\;\;8, \;8,\;8, \;8, \,9+)& 9\\
(3,\;\,5, 5, 6+) & 7 & (5,\;7, 7, 7,\;8+)& 9 &(7,\;\;9, \;9,\;9, \;9, 10+)& 11\\
(3,\;\,6, 6, 7+) & 7 & (5,\;8, 8, 8,\;9+)& 9 &(7,\;10, 10, 10, 10, 11+)& 11\\
(3,\;\.7, 7, 8+) & 9 & (5,\;9, 9, 9, 10+)&11 &(7,\;11, 11, 11, 11, 12+)& 13\\
\hline
\end{tabular}
\end{center}

\begin{center}
\footnotesize
\begin{tabular}{|ll|ll|}
\hline
$x_1 = 9, n = 7$ & $\R$ & $x_1 = 11, n = 8$ & $\R$\\
\hline
(9,\,\;\;9,\,\;\,9,\,\;\;\,9,\;\;\;9,\;\;\;\,9,\;\;10+) & 11 & (11,\;11, 11, 11, 11, 11, 11, 12+) & 13\\
(9,\,\;10, 10, 10, 10, 10, 11+) & 11 & (11,\;12, 12, 12, 12, 12, 12, 13+) & 13 \\
(9,\,\;11, 11, 11, 11, 11, 12+) & 13 & (11,\;13, 13, 13, 13, 13, 13, 14+) & 15 \\
(9,\,\;12, 12, 12, 12, 12, 13+) & 13 & (11,\;14, 14, 14, 14, 14, 14, 15+) & 15 \\
(9,\,\;13, 13, 13, 13, 13, 14+) & 15 & (11,\;15, 15, 15, 15, 15, 15, 16+) & 17\\
\hline
\end{tabular}
\end{center}


\subsection{Odd  $x_1 \geq 5$ with $n = \frac{1}{2}(x_1+1) + 1$}  
\label{ss16}  
Given an integer $i \geq 0$, 
a position $x = (x_1, \dots, x_n, \dots, x_m)$  is an exception 
if and only if one of the following two cases holds:  

$$x_1, \; x_2 = x_1 + 2i, \; x_3 = \dots = x_{n-1} = x_2 + 1 < x_n \leq \dots \leq x_m,$$  
$$x_1, \; x_2 = x_1 + 2i, \; x_3 = \dots = x_{n-1} = x_2 + 2 < x_n \leq \dots \leq x_m.$$  

Furthermore,  the remoteness function is given by formula 
$\R(x) = f(x_n)$, where function $f$ is defined above.  
Again, it is easily seen that  $\R(x)$  is odd for each  $i$; 
hence, the first player wins in every exceptional position, 
but (s)he loses if follows the M-rule. 

Finally,  
the unique optimal move in  $x$  is to keep  
$x_n$  when  $x_{n-1}$  is even and  
$x_{n-1}$   when it is odd. 
In contrast, the unique M-move in  $x$  is to keep  
$x_n$  when $x_{n-1}$ is odd and  
$x_{n-1}$  when it is even. 
Thus, the sets of optimal moves and M-moves 
are disjoint in every exceptional position.  

Note that  $x_3 = \ldots = x_{n-1}$  
and this number is even if and only if  $i$  is odd.

Examples for $x_1 = 5,7,9,11,13$  and, respectively, 
$n = 4,5,6,7,8$  are given below. 

\begin{center}
\footnotesize
\begin{tabular}{|ll|ll|ll|}
\hline
$x_1 = 5, n = 4$ & $\R$ & $x_1 = 7, n = 5$ & $\R$ & $x_1 = 9, n = 6$ & $\R$\\
\hline
(5, 5,\;\;6,\;7+) & \;7 & (7,\,7,\;\;8,\;\,8,\; 9+)& \, 9 & (9,\;9,\;\;10, 10, 10, 11+) & 11\\
(5, 5,\;\;7,\;8+) & \;9 & (7,\,7,\;\;9,\;\,9,\;10+)& 11 &(9,\;9,\;\;11, 11, 11, 12+) & 13\\
(5, 7,\;\;8,\;9+) & \;9 & (7,\,9,\;\,10, 10, 11+)& 11 &(9, 11,\;12, 12, 12, 13+) & 13\\
(5, 7,\;\;9, 10+) & 11 & (7,\,9,\,\;11, 11, 12+)& 13 &(9, 11,\;13, 13, 13, 14+) & 15\\
(5, 9,\;10, 11+) & 11 & (7, 11,\;12, 12, 13+)& 13 &(9, 13,\;14, 14, 14, 15+) & 15\\
\hline
\end{tabular}
\end{center}

\begin{center}
\footnotesize
\begin{tabular}{|ll|ll|ll|}
\hline
$x_1 = 11, n = 7$ & $\R$ & $x_1 = 13, n = 8$ & $\R$\\
\hline
(11, 11,\;\;\;12, 12, 12, 12, 13+) & 13 & (13, 13,\;\;\;14, 14, 14, 14, 14, 15+) &  15 \\
(11, 11,\;\;\;13, 13, 13, 13, 14+) & 15 & (13, 13,\;\;\;15, 15, 15, 15, 15, 16+) &  17 \\
(11, 13,\;\;\;14, 14, 14, 14, 15+) & 15 & (13, 15,\;\;\;16, 16, 16, 16, 16, 17+) & 17 \\
(11, 13,\;\;\;15, 15, 15, 15, 16+) & 17 & (13, 15,\;\;\;17, 17, 17, 17, 17, 18+) & 19 \\
(11, 15,\;\;\;16, 16, 16, 16, 17+) & 17 & (13, 17,\;\;\;18, 18, 18, 18, 18, 19+) & 19 \\
\hline
\end{tabular}
\end{center}

\subsection{Odd  $x_1 \geq 7$  with  $n = \frac{1}{2}(x_1+1)$}  
\label{ss17} 
Given an integer $i \geq 0$, 
a position $x = (x_1, \dots, x_n, \dots, x_m)$  is an exception 
if and only if one of the following two cases holds:  

$$x_1, \; x_2 = x_1 + 2i, \; x_3 = \ldots = x_{n-1} = x_2 + 3 < x_n \leq \ldots \leq x_m,$$  
$$x_1, \; x_2 = x_1 + 2i, \; x_3 = \ldots = x_{n-1} = x_2 + 4 < x_n \leq \ldots \leq x_m.$$  

Interestingly, all further arguments 
can be copied from the previous subsection 
without any changes, yet, 
we should remember that  $n$  is reduced by  1.  

Examples for $x_1 = 7, 9, 11, 13, 15$  and, respectively, 
$n = 4, 5, 6, 7, 8$  follow. 

\begin{center}
\footnotesize
\begin{tabular}{|ll|ll|ll|}
\hline
$x_1 = 7, n=4$ & $\R$ & $x_1 = 9, n=5$ & $\R$ & $x_1 = 11, n=6$ & $\R$\\
\hline
(7,\,7,\;\,10, 11+) & 11 & (9,\;9,\;12, 12, 13+) & 13&(11, 11,\;14, 14, 14, 15+) & 15\\
(7,\,7,\;\,11, 12+) & 13 & (9,\;9,\;13, 13, 14+) & 15&(11, 11,\;15,\,15, 15, 16+) & 17\\
(7,\,9,\;\,12, 13+) & 13 & (9,11,\;14, 14, 15+) & 15&(11, 13,\;16,\,16, 16, 17+) & 17\\
(7,\,9,\;\,13, 14+) & 15 & (9,11,\;15, 15, 16+) & 17&(11, 13,\;17,\,17, 17, 18+) & 19\\
(7,11,\;14, 15+) & 15 & (9,13,\;16, 16, 17+) & 17&(11, 15,\;18, 18, 18, 19+) & 19\\
\hline
\end{tabular}
\end{center}

\begin{center}
\footnotesize
\begin{tabular}{|ll|ll|}
\hline
$x_1 = 13, n=7$ & $\R$ & $x_1 = 15, n=8$ & $\R$\\
\hline
(13, 13,\;16, 16, 16, 16 ,17+) & 17 & (15, 15,\;18, 18, 18, 18, 18, 19+) & 19\\ 
(13, 13,\;17, 17, 17, 17, 18+) & 19 & (15, 15,\;19, 19, 19, 19, 19, 20+) & 21\\ 
(13, 15,\;18, 18, 18, 18, 19+) & 19 & (15, 17,\;20, 20, 20, 20, 20, 21+) & 21\\
(13, 15,\;19, 19, 19, 19, 20+) & 21 & (15, 17,\;21, 21, 21, 21, 21, 22+) & 23\\
(13, 17,\;20, 20, 20, 20, 21+) & 21 & (15, 19,\;22, 22, 22, 22, 22, 23+) & 23\\ 
\hline
\end{tabular}
\end{center}

There exists another family of exceptions 
for  $x_1 \geq 9$  with  $n = \frac{1}{2}(x_1+1)$. 
$$x_1, \; x_2 = x_3 = x_1 + 2i, \; x_4 = \ldots = x_{n-1} = x_2 + 3 < x_n \leq \ldots \leq x_m,$$  
$$x_1, \; x_2 = x_3 = x_1 + 2i, \; x_4 = \ldots = x_{n-1} = x_2 + 4 < x_n \leq \ldots \leq x_m.$$  
Examples for $x_1 = 9, 11, 13$  and, respectively, 
$n = 5, 6, 7$  follow. 

\begin{center}
\footnotesize
\begin{tabular}{|ll|ll|ll|}
\hline
$x_1 = 9, n=5$ & $\R$ & $x_1 = 11, n=6$ & $\R$ & $x_1 = 13, n=7$ & $\R$\\
\hline
(9,\;9, \;9,\; 10, 11+) & 11 & (11, 11, 11,\;12, 12, 13+) & 13 &(13, 13, 13,\;14, 14, 14, 15+) & 15\\
(9,\;9, \;9,\; 11, 12+) & 13 & (11, 11, 11,\;13, 13, 14+) & 15 &(13, 13, 13,\;15, 15, 15, 16+) & 17 \\
(9, 11, 11,\;12, 13+) & 13 & (11, 13, 13,\'14, 14, 15+) & 15 &(13, 15, 15,\;16, 16, 16, 17+) & 17 \\
(9, 11, 11,\;13, 14+) & 15 & (11, 13, 13,\;15, 15, 16+) & 17 &(13, 15, 15,\;17, 17, 17, 18+) & 19 \\
(9, 13, 13,\;14, 15+) & 15 & (11, 15, 15,\;16, 16, 17+) & 17 &(13, 17, 17,\;18, 18, 18, 19+) & 19\\
\hline
\end{tabular}
\end{center}

\subsection{Odd  $x_1 \geq 9$  with  $n = \frac{1}{2}(x_1 + 1) - 1$}  
\label{ss19} 


Exceptions for $x_1 = 9, 11, 13, 15$  and, respectively, 
$n = 4, 5, 6, 7$  follow: 

\begin{center}
\footnotesize
\begin{tabular}{|ll|ll|ll|}
\hline
$x_1 = 9, n = 4$ & $\R$ & $x_1 = 11, n = 5$ & $\R$ & $x_1 = 13, n = 6$ & $\R$ \\
\hline
(9, \; 9, \; 14, 15+) & 15 &(11, 11, 13,\;14, 15+) &  15 & (13, 13, 13, 13,\;14, 15+) & 15  \\
(9, \; 9, \; 15, 16+) & 17 &(11, 11, 13,\;15, 16+) &  17 & (13, 13, 13, 13,\;15, 16+) & 17 \\
(9, 11, \; 16, 17+) & 17 &(11, 13, 15, \;16, 17+) &  17 & (13, 15, 15, 15,\;16, 17+) & 17 \\
(9, 11, \; 17, 18+) & 19 &(11, 13, 15, \;17, 18+) &  19 & (13, 15, 15, 15,\;17, 18+) & 19 \\
(9, 13, \; 18, 19+) & 19 & (11, 15, 17, \;18, 19+) &  19 & (13, 17, 17, 17,\;18, 19+) & 19  \\ 
\hline
\end{tabular}
\end{center}

\begin{center}
\footnotesize
\begin{tabular}{|ll|ll|}
\hline
$x_1 = 13, n = 6$ & $\R$ & $x_1 = 15, n = 7$ & $\R$\\
\hline
(13, 13, 15, \; 16, 16, 17+) & 17& (15, 15, 15, 15, \; 16, 16, 17+) & 17 \\ 
(13, 13, 15, \; 17, 17, 18+) & 19& (15, 15, 15, 15, \; 17, 17, 18+) & 19 \\
(13, 15, 17, \; 18, 18, 19+) & 19 & (15, 17, 17, 17, \; 18, 18, 19+) & 19\\
(13, 15, 17, \; 19, 19, 20+) & 21& (15, 17, 17, 17, \; 19, 19, 20+) & 21 \\
(13, 17, 19, \; 20, 20, 21+) & 21& (15, 19, 19, 19, \; 20, 20, 21+) & 21 \\
\hline
\end{tabular}
\end{center}

\medskip 

\subsection{Odd  $x_1 \geq 11$  with  $n = \frac{1}{2}(x_1 + 1) - 2$}  
\label{ss19a} 

For  $n = (x_1 + 1)/2 - 2$   
we obtained the following exceptions: 
\bigskip 

\begin{center}
\footnotesize
\begin{tabular}{|ll|ll|l|}
\hline
$x_1 = 11, n=4$ & $\R$ & $x_1 = 13, n=5$ & $\R$ & $x_1 = 15, n=6$\\
\hline
(11, 11, \; 18, 19+) & 19&(13, 13, \; 20, 20, 21+) & 21&(15, 15, \; 22, 22, 22, 23+)\\
(11, 11, \; 19, 20+) & 21&(13, 13, \; 21, 21, 22+) & 23&(15, 17, \; 23, 23, 23, 24+)\\
(11, 13, \; 20, 21+) & 21&(13, 15, \; 22, 22, 23+) & 23&(15, 17, \; 24, 24, 24, 25+)\\ 
(11, 13, \; 21, 22+) & 23&(13, 15, \; 23, 23, 24+) & 25&(15, 19, \; 25, 25, 25, 26+)\\
(11, 15, \; 22, 23+) & 23&(13, 17, \; 24, 24, 25+) & 25&(15, 19, \; 26, 26, 26, 27+)\\
\hline
\end{tabular}
\end{center}

In addition, the following exceptions were found: 

\begin{table}[h]
\footnotesize
\begin{tabular}{|ll|}
\hline
$x_1 = 13, n=5$ & $\R$ \\
\hline
(13, 13, 13,\; 16, 17+) & 17\\
(13, 13, 13,\; 17, 18+) & 19\\
(13, 15, 15,\; 18, 19+) & 19\\
(13, 15, 15,\; 19, 20+) & 21\\
(13, 17, 17,\; 20, 21+) & 21\\
\hline
\end{tabular}
\hfill
\begin{tabular}{|p{3.2cm}|l|}
\hline
$x_1 = 15, n=6$\\
\hline
(15, 15, \; 23, 23, 23, 24+)\\
(15, 17, \; 24, 24, 24, 25+)\\
(15, 17, \; 25, 25, 25, 26+)\\
(15, 19, \; 26, 26, 26, 27+)\\
(15, 19, \; 27, 27, 27, 28+)\\
\hline
\end{tabular}                   
\hfill
\begin{tabular}{|p{3.9cm}|l|}
\hline
$x_1 = 17, n=7$\\
\hline
(17, 17, 17, 17, 17, \; 18, 19+)\\
(17, 17, 17, 17, 17, \; 19, 20+)\\
\\
\\
\\
\hline
\end{tabular}                
\end{table}

\subsection{Odd  $x_1 \geq 13$  with  $n = \frac{1}{2}(x_1 + 1) - 3$}  
\label{ss19a} 

Exceptions for $x_1 = 13, 15, 17$  and, respectively, 
$n = 4, 5, 6$  follow: 

\begin{center}
\footnotesize
\begin{tabular}{|l|l|l|}
\hline
$x_1 = 13, n=4$ & $x_1 = 15, n=5$ & $x_1 = 17, n=6$ \\
\hline
(13, 13, \; 22, 23+)&(15, 15, \; 24, 24, 25+)&(17, 17, \; 26, 26, 26, 27+)\\ 
(13, 13, \; 23, 24+)&(15, 15, \; 25, 25, 26+)&(17, 17, \; 27, 27, 27, 28+) \\
(13, 15, \; 24, 25+)&(15, 17, \; 26, 26, 27+)&(17, 19, \; 28, 28, 28, 29+) \\
(13, 15, \; 25, 26+)&(15, 17, \; 27, 27, 28+)&(17, 19, \; 29, 29, 29, 30+) \\
(13, 17, \; 26, 27+)&(15, 19,\; 28, 28, 29+)&(17, 21, \; 30,  30, 30, 31+)\\ 
\hline
\end{tabular}
\end{center}

\bigskip 
\noindent 
{\bf Acknowledgements.}
This research was prepared within the framework 
of the HSE University Basic Research Program.  

\end{document}